\documentclass[a4paper,11pt,draft]{article}
\usepackage{graphics}
\usepackage{amsmath}
\usepackage{amsfonts}
\usepackage{amssymb}
\usepackage[russian]{babel}
\usepackage{amscd}
\newtheorem{theorem}{\sc \bf Theorem.}
\newenvironment{proof1} {\smallskip\normalfont{\scshape \bf Proof of the theorem 1.}}{\vspace{.2 ex}}
\newenvironment{proof2} {\smallskip\normalfont{\scshape \bf Proof of the theorem 2.}}{\vspace{.2 ex}}
\newenvironment{proof}{{\bf Proof.}}{\par\hspace{25em}\rule{1eX}{1eX}\par}
\newtheorem{lemma}{Лемма.}
\newcommand{\vr}{\varphi}
\newcommand{\opred}[1]{\it #1}
 \newcommand{\ds}{\displaystyle}
\begin{document}
\begin{center}
 {\bf Convex hypersurfaces in Hadamard manifolds.}
 \end{center}
\begin{center}
 {\bf A.A.Borisenko.}
   \end{center}

\hrule \vspace{0.5cm} {\bf Abstract.} We prove the theorem about
extremal property of Lo\-ba\-chev\-sky space among simply
connected Riemannian manifolds of non\-po\-si\-ti\-ve
cur\-va\-ture. \vspace{0.3cm}
 \hrule
 \vspace{0.5cm}
 Hadamard proved the following theorem. Let $\vr$ be an immersion of
 a com\-pact ori\-ented $n$-dimensional manifold $M$ in Euclidean space
 $E^{n+1} \, (n\geqslant~2)$ with everywhere positive Gaussian
 curvature. Then $\vr(M)$ is a convex hy\-per\-sur\-fa\-ce~[1].

 Chern and Lashof [2] generalized this theorem. Let $\vr$
 be an immersion of a compact oriented $n$-dimensional manifold $M$
 in $E^{n+1}$. Then the following two assertions are equivalent:
\begin{enumerate}
\item[(i)] The degree of the spherical mapping equals $\pm1$, and
the Gaussian curvature does not change sign (i.e., it is
everywhere nonnegative or everywhere nonpositive); \item[(ii)]
 $\vr(M)$ is a convex hypersurface.
\end{enumerate}
By Gaussian curvature, we mean the product of the principal
curvatures.

S.Alexander generalized Hadamard theorem for compact hypersurfaces
in any complete, simply connected Riemannian manifold of
non\-po\-si\-ti\-ve sec\-tio\-nal cur\-va\-ture.\cite{3}

A topological immersion $f: N^n\to M$ of a manifold $N^n$ into a
Rie\-ma\-nni\-an ma\-ni\-fold $M$ is called {\opred locally
convex} at a point $x\in N^n$ if  has a neighbourhood $U$ such
that $f(U)$ is a part of the boundary of a convex set in $M$.

Heijenoort proved the following theorem. Let $f: N^n\to E^{n+1}$,
where $n\geq2$, be a topological immersion of a connected manifold
$N^n$. If $f$ is locally convex at all points and has at least one
point of local strict support and $N^n$ is complete in the metric
induced by immersion, then $f$ is an embedding and $F=f(N^n)$ is
the boundary of a convex body \cite{4}.

In \cite{5}   this theorem was generalized to $h$-locally convex
(i.e., such that each point has a neighbourhood lying on one side
from a horosphere) regular hypersurfaces in Lobachevsky space and
in \cite{6}, to nonregular hypersurfaces.

In this section we shall recall some definitions and we shall
state the notation.

A {\opred Hadamard manifold} is a complete simply connected
Riemannian manifold  with sectional curvature $K\leq0$.

Like in the hyperbolic space, a {\opred horoball in a Hadamard
manifold } $M$ is the domain obtained as the limit of the balls
with their centres in a geodesic ray going to infinity, and their
corresponding geodesic spheres containing a fixed point. The
boundary of a horoball is a {\opred horosphere}. In general, a
horosphere is a $C^2$ hypersurface. An {\opred $h$-convex set} in
a Hadamard manifold $M$ of dimension $n+1$ is a subset
$\Omega\subset M$ with boundary $\partial\Omega$ satisfying that,
for every $P\in\partial\Omega$ there is a horosphere $H$ of $M$
through $P$ such that $\Omega$ is locally contained in the
horoball of $M$ bounded by $H$. This $H$ is called a {\opred
supporting horosphere} of $\Omega$ (and $\partial\Omega$).

For Hadamard manifolds $M$ satisfying $-k_1^2\geqslant
K\geqslant-k_2^2,\, k_1, k_2>0$, if $H$  is horosphere, at each
point of $H$ where the normal curvature $k_n$ is well defined, it
satisfies $k_1\leq k_n\leq k_2$.

For geodesic spheres of radius $r$ normal curvatures satisfy
inequality
$$k_1\coth k_1r\leq k_n\leq k_2\coth k_2r.$$
Note that the value $k\coth kr$ is the geodesic curvature of a
circumference of radius $r$ in Lobachevsky plane of curvature
$-k^2$.

{\opred An orientable regular ($C^2$ or more) hypersurface $F$ of
a Hadamard manifold $M$ is $\lambda$-convex if, for a selection of
its unit normal vector, the normal curvature $k_n$ of $F$
satisfies $k_n\geqslant\lambda$. A domain $\Omega\subset M$ is
$\lambda$-convex if for every point $P\in\partial\Omega $ there is
a regular $\lambda$-convex  hypersurface $F$ through $P$ leaving a
neighbourhood of~$P$  in the convex side (the side where the unit
normal vectors points) of $F$.} If $\partial\Omega$ is regular,
then it is a regular $\lambda$-convex hypersurface.

Given any set $\Omega\subset M $, an {\opred inscribed ball
(inball} for short) is a ball in $M$ contained in $\Omega$ with
maximum radius. Its radius is called the {\opred inradius} of
$\Omega$, and it will be always denoted by $r$. Moreover, we shall
denote by $O$ the (not necessarily unique) centre of an inball of
$\Omega$, and by $d$ the distance, in $M$ to~$O$.

A circumscribed ball (or circumball) is a ball in $M^{n+1}$
containing $\Omega$ with minimum radius. Its radius is called
circumradius  in $\Omega$ and notes by $R$.

Now we shall prove the following theorems.
\begin{theorem}\label{1}
Let $M^{n+1}$ be a simply connected complete Riemannian manifold
with sectional curvature
  $$-k_1^2\geqslant K\geqslant -k_2^2,\quad k_2\geqslant k_1>0.$$

  Suppose that $F\subset M^{n+1}$  be a complete immersed
  hypersurface with nor\-mal curvatures
  $$k_n\geqslant k_2.$$
  Then either
  \begin{enumerate}
  \item[I)] $F^n$ is a compact convex hypersurface diffeomorphic to
  the sphere $S^n$ and
  $$R-r<k_2\ln2 $$
  or
\item[II)] $F^n$ is a horosphere in $M^{n+1}$ and ambient space
$M^{n+1}$ is a hyperbolic space of constant curvature $-k_2^2$.

For more strong condition on the normal curvatures $F^n$ it is
true
\end{enumerate}
\end{theorem}
\begin{theorem}
Let $M^{n+1}$ be a Hadamard manifold with sectional curvature
  $$-k_1^2\geqslant K\geqslant -k_2^2,\quad k_2\geqslant k_1\geqslant0.$$
  Let $F^n$ be a complete immersed hypersurface with normal
  curvatures more or equal
$k_2\coth k_2r_0$ at any point $F^n$. Then
  either
  \begin{enumerate}
  \item[I)] $F^n$ is a compact convex hypersurface diffeomorphic
  sphere $S^n$ and radius of circumscribed ball
  $$R<r_0$$
  or
  \item[II)]
   $F^n$ is a sphere of radius $r_0$ which is boundary the ball
   $\Omega$. The ball $\Omega$ is isometric to the ball of radius
   $r_0$ of hyperbolic space with constant curvature $-k_2^2$
  \end{enumerate}
\end{theorem}

 The ambient space $M^{n+1}$ is a $C^3$ regular Riemannian
 manifold. For proof of part $II)$ of theorem \eqref{1} we need
 the condition.
$$|\nabla R|\leq C,$$
where $|\nabla R|$ is a covariant differential of curvature tensor
$M^{n+1}$, $C$ is a positive constant.

For condition $K_{\sigma}\leq -k^2_1<0$ and $|\nabla R|< C$ a
horosphere is $C^3$ re\-gu\-lar hy\-per\-sur\-fa\-ce and
ma\-ni\-fold $M^{n+1}$ is $C^2$-regular Riemannian  manifold in
horospheric coordinates \cite{7}.

At any point smooth hypersurface $F^n$ in Hadamard manifold there
are two tangent horospheres. Let normal curvatures $F^n$  at some
point $P\in F^n$ with respect some normal be greater zero, one of
the horospheres with positive normal curvature with respect the
same normal we call {\opred tangent horosphere.}

\begin{proof1}
\begin{enumerate}
\item[I).] From the condition of the theorem it follows that
normal curvature of the horosphere $H^n$ in $M^{n+1}$
$$ k_n/_{H^n}\leq k_2.$$ And for every point $P\in F^n$, normal
curvatures of  tangent horosphere   in \newpage the  corresponding
directions satisfy the inequality
$$k_n(a)/_{F^n}\geqslant k_n(a)/_{H^n}.$$
\quad Suppose that in point $P_0$ it is true the strong inequality
\begin{equation}\label{1.1}
k_n(a)/_{F^n}> k_n(a).
\end{equation}
\quad Let $n_0$ be the unit normal at the point $P_0$, such that
the normal curvatures of $F^n$ at the point $P_0\in F^n$ with
respect normal $n_0$ are positive, $H^n$ be a  tangent horosphere
at the point $P_0$ with the normal $n_0$. From the
inequality~\eqref{1.1} it follows that there exists some
neighbourhood of the point $P_0$ on $F^n$ such that it lies inside
the horoball bounded by horosphere $H^n$. Let we take a
horospherical system of the coordinates in $M^{n+1}$ with the base
$H^n$, $t$ is a length parameter along geodesic line orthogonal
$H^n$, positive direction coincides with the normal  $n_0$ at the
point $P_0$. From another side $t$ is a distance  from a point of
$M^{n+1}$ to the horosphere $H^n$. Let the function $f=t$ be the
restriction $t$ oh the hypersurface $F^n$, at the point $P_0$ the
function $f$ has a strong minimum. Let $\vr$ be the angle  between
the direction $\frac{\partial}{\partial t}$ and the unit normal
$N$ of the hypersurface $F^n$. Along integral curves of the vector
field $X=grad\, f/_{F^n}$ on the hypersurface $F^n$ the angle
$\vr$ satisfies the equation \cite{8}.
\begin{equation}\label{2.2}
k_n=\mu\cos\varphi+\sin\varphi\frac{d\varphi}{dt},
\end{equation}
where $k_n$ is the normal curvature $F^n$ in the direction $X$ at
the point $P\in F^n$, $\mu$ is the normal curvature of the
coordinate horosphere at the point $P\in F^n$ in the direction
$Y$, which is orthogonal projection the vector $X$ on the tangent
space of the coordinate horosphere at the point $P$.

As $k_n\geqslant k_2$ and normal curvatures of the horosphere
$\mu\leq k_2$ that from \eqref{2.2} for $\vr\leq\frac{\pi}{2}$ it
follows
$$
\begin{array}{c}
k_2(1-\cos\varphi)\leq\sin\varphi\frac{d\varphi}{dt};\\[2ex]
\dfrac{\sin\frac{\varphi}{2}}{\sin\frac{\varphi_0}{2}}\geqslant e^{\frac{k_2}{2}(t-t_0)};\\[2ex]
\sin\frac{\varphi}{2}\geqslant\sin\frac{\varphi_0}{2}
e^{\frac{k_2}{2}(t-t_0)},
\end{array}
$$
where $\vr_0>0$ is the angle between $\frac{\partial}{\partial t}$
and the normal $N$ for small $t_0$. It follows from inequality
\eqref{1.1} at the point $P_0$. The angle $\vr$ monotonically
increases along integral curve and for
$$t\leq\frac{2}{k_2}\ln{\frac{e^{\frac{k_2t_0}{2}}}{\sqrt{2}(\sin\frac{\varphi_0}{2})}}$$
reaches the value $\frac{\pi}{2}$. For $\vr\geqslant\frac{\pi}{2}$
we have
$$
\begin{array}{c}
k_2\leq\sin\varphi\dfrac{d\varphi}{dt};\\[2ex]
\cos\varphi\leq 1-k_2(t-t_1),
\end{array}
$$
where $\vr(t_1)=\frac{\pi}{2}$ and for $t_2\leq t_1+\frac{2}{k_2}$
the angle $\vr$ reaches the value $\pi$ and function $f=t/_{F^n}$
at this point achieves strong maximum.

The length of integral curve on the hypersurface $F^n$ of the
vector field $X~=~grad\, f/_{F^n}$ satisfies the inequality
$$
\begin{array}{c}
k_2(1-\cos\varphi)\leq\dfrac{d\varphi}{ds};\\[2ex]
s\leq s_0+\dfrac{\cot{\frac{\vr_0}{2}}}{k_2}.
\end{array}
$$
It follows that point $Q_0$, where $\vr=\pi$ does not go to
infinity. Let $t_2$ be the infinum of the value $t$ on integral
curves of vector field $X=grad\, f/_{F^n}$ such that
$\vr(t_2)=\pi$. Level hypersurfaces of the fun\-cti\-on $f=t$ for
$0<t<t_2$ are sphe\-res $S^{n-1}$ and points $P_0$ and $Q_0$ are
strong minimum and strong maximum and function $f$ is a Morse
function on $F^n$ with two critical points. Therefore the
hypersurface $F^n$  is homeomorphic to sphere $S^n$. From the
condition $k_n~\geqslant~k_2$ we obtain that second quadratic form
$F^n$ is po\-si\-ti\-ve de\-fi\-ni\-te at any point. From theorem
S. Alexander \cite{4} it follows that $F^n$ is embedded compact
convex hypersurface diffeomorphic  $S^n$ and bounds convex domain
$\Omega$. The domain $\Omega$ is $k_2$-convex and satisfies the
condition of theorem $3.1$~\cite{9} and
$$\max\, d(O, \partial\Omega)<k_2\ln2,$$
where $O$ is the centre of the inscribed ball.
\begin{enumerate}
\item[II)1).] Suppose that at any point $P\in F^n$  there exists
the direction
 $a~\in~T_pF^n$ such that
 $$k_n(a)/_{F^n}=k_n(a)/_{H^n}.$$
Let show that some neighbourhood $U\subset F^n$ of a point $P_0\in
 F^n$ lies in the horoball bounded by tangent horosphere $H^n$. Let
 take horospherical system of coordinate with the base $H^n$.

 The metric $M^{n+1}$ has the form
 \begin{equation}\label{bor2}
ds^2=dt^2+g_{ij}(t,\theta)d\theta^id\theta^j.
\end{equation}
The equation of the hypersurface $F^n$ in the neighbourhood
$P_0\in F^n$ is
$$t=\rho(\theta).$$
The unit normal vector $N$ to $F^n$ has coordinates
\begin{equation}
\xi^k=-\dfrac{\rho^k}{\sqrt{1+\langle
grad\,\rho,grad\,\rho\rangle}},\quad {k=1,\ldots,n}
\end{equation}
$$
\xi^{n+1}=\dfrac{1}{\sqrt{1+\langle
grad\,\rho,grad\,\rho\rangle}},
$$
where
\begin{equation}
\rho^k=g^{ks}\rho_s,\quad \langle
grad\,\rho,grad\,\rho\rangle=g_{ij}\rho^i\rho^j,
\end{equation}
$$\rho_i=\frac{\partial\rho}{\partial\theta^i}.$$
Coefficient of the second fundamental form of $F^n$  is equal
\cite{10}
\begin{equation}\label{B5}
\Omega_{ij}=\cos\varphi\left[\rho_{i,j}-\dfrac{1}{2}\dfrac{\partial
g_{ij}}{\partial t}-\dfrac{1}{2}\dfrac{\partial g_{jk}}{\partial
t}\rho_i\rho^k-\dfrac{1}{2}\dfrac{\partial g_{ik}}{\partial
t}\rho_j\rho^k\right],
\end{equation}
where $\vr$ is the angle between $\frac{\partial}{\partial t}$ and
normal $N$,
\begin{equation}
\cos\varphi=\dfrac{1}{\sqrt{1+\langle
grad\,\rho,grad\,\rho\rangle}},\quad
\rho_{i,j}=\rho_{ij}-\Gamma^k_{ij/g}\rho_k,
\end{equation}
 where $\Gamma^k_{ij/g}$ are Kristoffel symbols of the
metric
$$d\sigma^2=g_{ij}d\theta^i\,d\theta^j.$$

Coefficients of metric tensor $F^n$ have the form
\begin{equation}
a_{ij}=g_{ij}+\rho_i\rho_j.
\end{equation}

From the conditions of the theorem normal curvatures
$$k_n/_{F^n}\geqslant k_2.$$

And from \eqref{B5} it follows that for any tangent vector $b\in
F^n,\,\\ b=(b^1,\ldots, b^n).$
\begin{multline}\label{B8}
\cos\varphi\left[\rho_{i,j}b^ib^j-\dfrac{1}{2}\dfrac{\partial
g_{ij}}{\partial t}b^ib^j-\dfrac{1}{2}\dfrac{\partial
g_{jk}}{\partial t}\rho_i b^i\rho^kb^j-\dfrac{1}{2}\dfrac{\partial
g_{ik}}{\partial t}\rho_jb^j\rho^kb^i\right]\geqslant\\\geqslant
k_2(g_{ij}b^ib^j+(\rho_ib^i)^2).
\end{multline}
Let introduce the function $h=e^{k_2\rho(\theta)}$.
$$
\begin{array}{c}
h_i=k_2e^{k_2\rho}\rho_i;\\[2ex]
h_{ij}=k_2^2e^{k_2\rho}\rho_i\rho_j+k_2e^{k_2\rho}\rho_{ij}.
\end{array}
$$
Hence
$$\rho_i=\dfrac{h_i}{h}\dfrac{1}{k_2};$$
\begin{equation}\label{B9}
\rho_{ij}=\dfrac{1}{k_2}\dfrac{h_{ij}}{h}-\dfrac{1}{k_2}\dfrac{h_i}{h}\dfrac{h_j}{h};
\end{equation}
$$\rho_{i,j}=\dfrac{1}{k_2}\dfrac{hh_{i,j}-h_ih_j}{h^2}.$$
And inequality \eqref{B8} we rewrite in the following way:
\begin{equation}\label{11}
\begin{array}{r}
\ds \cos\varphi\left[
\dfrac{1}{k_2}hh_{i,j}b^ib^j-\dfrac{1}{k_2}(h_ib^i)^2
-\dfrac{1}{2}\dfrac{\partial g_{jk}}{\partial t}b^i
b^jh^2-\dfrac{1}{2k_2^2}\dfrac{\partial g_{ik}}{\partial
t}h_ib^ih^kb^j- \right.\\[3ex]
\ds \left.\dfrac{1}{2k_2^2}\dfrac{\partial g_{ik}}{\partial
t}h_jb^jh^kb^i\right]\geqslant
k_2\left[h^2g_{ij}b^ib^j+\dfrac{1}{k_2^2}(h_ib^i)^2\right]
\end{array}
\end{equation}
Since the normal curvature of horosphere in $M^{n+1}$ less or
equal $k_2$ than
\begin{equation}\label{B12}
-\dfrac{1}{2}\dfrac{\partial g_{ij}}{\partial t}b^ib^j\leq
k_2g_{ij}b^ib^j,
\end{equation}
where $A_{ij}=-\dfrac{1}{2}\dfrac{\partial g_{ij}}{\partial t}$
are coefficients of the second fundamental form of horosphere
$t=const$.
\begin{multline}\label{Bor14}
\left|-\dfrac{1}{2}\dfrac{\partial g_{jk}}{\partial t}h
_ib^ih^kb^j\right|=(A_{jk}h^kb^j)|h_ib^i|\leq\\
\leq\sqrt{(A_{jk}h^kh^j)A_{jk}b^kb^j}|h_ib^i|\leq k_2|grad\,
h||b||h_ib^i|,
\end{multline}
where $|b|^2=g_{ij}b^ib^j$, \quad $|grad\,h|^2=g_{ij}h^ih^j.$

Let we substitute \eqref{B12}, \eqref{Bor14} in \eqref{11} and
obtain
\begin{multline}\label{B15}
\cos\varphi\dfrac{1}{k_2}hh_{i,j}b^ib^j\geqslant
k_2h^2(1-\cos\varphi)|b|^2+\dfrac{1}{k_2}(1+\cos\varphi)(h_ib^i)^2-\\
-2\dfrac{1}{k_2}|grad\,h||b||(h_ib^i)|.
\end{multline}
The expression in the right side in the quadratic equation with
respect $|(h_ib^i)|$. The discriminant of this equation is
\begin{equation}\label{B16}
\dfrac{1}{k_2^2}|grad\,h|^2|b|^2-h^2\sin^2\varphi|b|^2.
\end{equation}
But $$
\cos^2\varphi=\dfrac{1}{1+|grad\,\rho|^2}=\dfrac{k_2^2h^2}{k_2^2h^2+|grad\,h|^2},
$$
$$
\sin^2\varphi=\dfrac{|grad\,h|^2}{k_2^2h^2+|grad\,h|^2}.$$ And we
rewrite \eqref{B16} in the form
\begin{equation}\label{B18}
\dfrac{|b|^2}{k_2^2}\left(\dfrac{|grad\,h|^4}{k_2^2h^2+|grad\,h|^2}\right)\geqslant0.
\end{equation}
From \eqref{B15} it follows
 \begin{equation}\label{B19}
h_{i,j}b^ib^j\geqslant0.
\end{equation}
Let $L$ be lines on $F^n$ which satisfy the system of the
equations
\begin{equation}\label{B20}
\dfrac{\partial^2\theta^k}{\partial
s^2}+\Gamma^k_{ij/g}(\theta,\rho(\theta))\dfrac{\partial\theta^i}{\partial
s}\dfrac{\partial\theta^s}{\partial s}=0.
\end{equation}
From any point and in any direction goes through only one line
from this family. These line we call g-geodesic. We take the
restriction the function $h$ on this line $$
\begin{array}{c}
\theta^i=\theta^i(s);\\
h_s=h_i\dfrac{d\theta^i}{ds};
\end{array}
$$
\begin{equation}\label{B21}
h_{ss}=h_{ij}\dfrac{d\theta^i}{ds}\dfrac{d\theta^j}{ds}+h_k\dfrac{d^2\theta^k}{d^2s}.
\end{equation}
If we substitute \eqref{B20} in \eqref{B21} then
\begin{equation}\label{B22}
h_{ss}=h_{i,j}\dfrac{d\theta^i}{ds}\dfrac{d\theta^j}{ds}\geqslant0.
\end{equation}
At the point $P_0\quad h=1, h_s=0$ and from \eqref{B22} it follows
that along g-geodesic lines which go through the point $P_0,\,
h\geqslant1$. On tangent horosphere $H^n, \, h=1$ and the
hypersurface $F^n$ lies from one side tangent horosphere $H^n$.
 \item[2).]\quad Let $P_0$ be an arbitrary  fixed point $F^n$,
 $H^n(P_0)$~---~tangent horosphere. From $1)$ it
 follows that some neighbourhood of the point $P_0\in F^n$ is
 situated in horoball bounded by horosphere $H^n(P_0)$. Let take
 dual tangent horosphere $\tilde{H}^n(P_0)$. This horosphere is
 defined by opposite point at infinity on geodesic line going in
 the direction of normal $n_0$ at the point $P_0\in
 F^n,\,\tilde{H}^n(\tau)$ are parallel horospheres $\tilde{H}(0)~=~\tilde{H}^n(P_0),\,
  M_{\tau}=F^n\bigcap\tilde{H}^n(\tau)$, $\tau$ is a distance from the
  horosphere $\tilde{H}^n(P_0)$. Let $f=\tau/_{F^n}$ is the
  restriction the function $\tau$ on the hypersurface $F^n$. For
  the function $f$ the point $P_0$ is a strong local minimum for
  small $\tau$ the set $ M_{\tau}=F^n\bigcap\tilde{H}^n(\tau)$ is
  a diffeomorphic to the sphere $S^{n-1}$ and bounds on $F^n$ the
  domain $D_{\tau}$ homeomorphic a ball and contains unique
  critical point $P_0$ of the function  $f=\tau/_{F^n}$. On the
  horosphere $\tilde{H}^n(\tau)$ the set  $M_{\tau}$ bounds convex
  domain homeomorphic a ball. Really, the normal $\nu$ to
  $M_{\tau}$ on $\tilde{H}^n(\tau)$ has the form
  $$\nu=\lambda_1n_1+\lambda_2N,$$
  where $n_1$ is unit normal to $\tilde{H}^n(\tau)$, $N$ is a normal to
  $F^n, \langle\nu, n_1\rangle=0$. Therefore
  $$\nu=\langle n_1,N\rangle n_1+N.$$
  Let $X$ be the unit vector field tangent to $M_{\tau}$. Then
  $$\langle \nu,\nabla_XX\rangle=\langle n_1,N\rangle\mu+k_n/_{F^n},$$
  where $\mu$ is the normal curvature of the horosphere $\tilde{H}(\tau)$.
Since ${k_n/_{F^n}}~{\geqslant~k_2}$ and $\mu\leq k_2$,then
$\langle \nu,\nabla_XX\rangle>0$,that is the second quadratic form
$M_\tau$ on $\tilde{H}^n(\tau)$ is a positive definite and the
domain on $\tilde{H}^n(\tau)$ bounded $M_\tau$  is a convex domain
homeomorphic a ball.

\quad Let's consider the body  $Q(\tau)$, bounded
$\mathcal{D}_\tau$ and $\tilde{H}^n(\tau)$ for small $\tau$. At
any boundary point there  exists  a local  supporting horosphere.
It is a global supporting horosphere too.  And the body $Q(\tau)$
is situated in the horoball bounded by supporting horospheres.
Other words the body $Q(\tau)$ is
 $h$-convex. Let $\tau^*$ be a supremum $\tau$, for which the body $Q(\tau)$
is $h$-convex, $\mathcal{D}^*=\bigcup\mathcal{D}_\tau$. Let's show
$\tau^*=\infty$. Let us assume the contrary. There are three
possible cases:
\begin{enumerate}
\item[a).] $\mathcal{D}^*=F^n$; \item[b).] $\mathcal{D}^*\neq F^n$
and on the boundary $S^*$ of the domain $\mathcal{D}^*$ there are
critical points the function $f=\tau/F^n$. \item[c).]
$\mathcal{D}^*\neq F^n$ and  $S^*$ doesn't contain the critical
points  the function~$f$.
\end{enumerate}

\quad The case $c)$ is impossible. Really for $\tau>\tau^*$ the
set $M_\tau$ is ho\-me\-omor\-phic the sphere too. It bounds the
convex domain on  $\tilde{H}^n(\tau)$ and at any boundary point
$Q(\tau)$ there exists  a local supporting  horosphere. It follows
that $Q(\tau)$ is a $h$-convex set for $\tau>\tau^*$ and $\tau^*$
is not  supremum.

\quad At the case $b)$ the set $S^*$ contains a critical point $P$
of function $f$. At point $P\in S^*$ the horosphere
$\tilde{H}^n(\tau^*)$ is the tangent supporting horosphere to
$F^n, S^*\subset \tilde{H}^n(\tau^*)\bigcap F^n$ is the boundary
of the convex domain homeomorphic a ball on $\tilde{H}^n(\tau^*)$.
Let show that $\tilde{H}^n(\tau^*)$ is the tangent horosphere at
all points $S^*$. Really, some neighbourhood  $U$ of the point
$P\in F^n$ lies at one side with respect to $\tilde{H}^n(\tau^*),
U\bigcap S^*$  belongs $\tilde{H}^n(\tau^*)$. If the horosphere
$\tilde{H}^n(\tau^*)$ isn't tangent in some point  $Q\in U\bigcap
S^*$ then $U$ doesn't lies for one side $\tilde{H}^n(\tau^*)$. The
set $S^*$ is homeomorphic to the sphere $S^{n-1}$ and the sets of
the points of $S^*$, such that the horosphere
$\tilde{H}^n(\tau^*)$ is tangent, is opened and closed at the same
time. This set isn't empty and coincides with $S^*$. Let
$Q(\tau^*)$ be the body bounded $\mathcal{D}^*$  and the domain
with boundary $S^*$ on $\tilde{H}^n(\tau^*)$. It is a compact
$h$-convex body with smooth boundary. Let $S(r)$ be the
circumscribed sphere $Q(\tau^*)$ .

\quad Suppose that a tangent point $P\in S(r)$ to the boundary
$Q(\tau^*)$ belongs to $\tilde{H}^n(\tau^*)$. At this case the
sphere $S(r)$ is supporting to the horosphere
$\tilde{H}^n(\tau^*)$ at the point $P$. The sphere $S(r)$ and
$\tilde{H}^n(\tau^*)$ are tangent at the point $P$ and convex
sides have the same direction. This is impossible.

\quad Let $Q_0\in\mathcal{D}^*$ be a tangent point of the sphere
$S(r)$. For Hadamard manifolds are true the following.
\begin{lemma}
Let  $S(r)$ and  $S(R)\,(r<R)$be tangent spheres at the point $Q$
in Hadamard manifold of the sectional curvatures $K\leq0$.

\quad Suppose at the point $Q$ the convex sides of the spheres
are the same. Then at the point $Q$ the normal curvatures of the
sphere $S(R)$are less than normal curvatures the sphere $S(r)$ in
corresponding directions.
\end{lemma}

\begin{proof}
Let take in   $M^{n+1}$ the spherical system coordinate with pole
$O$, where $O$ is the centre  of the sphere  $S(R)$. In the
neighbourhood of the point $Q$ the sphere $S(r)$ has the following
parametrization
$$t=h(\theta^1,\ldots,\theta^n),$$
where $t,\theta^1,\ldots,\theta^n$ are spherical coordinates in
$M^{n+1}$ with metric
$$ds^2=dt^2+g_{ij}(t,\theta)d\theta^id\theta^j.$$
The normal curvature $S(r)$ at the tangent point $Q$ of the
spheres in the direction  $b=(b^1,\ldots,b^n)$ is equal:
\begin{multline}\label{Bor19}
k_n=\dfrac{\left(\dfrac{\partial^2h(\theta^1,\ldots,\theta^n)}{\partial
\theta^i\partial\theta^j}-\dfrac{1}{2}\dfrac{\partial
g_{ij}}{\partial
t}\right)b^ib^j}{g_{ij}b^ib^j}=\\
k_n(b)/_{S(R)}+\dfrac{\dfrac{\partial^2h(\theta^1,\ldots,\theta^n)}
{\partial\theta^i\partial\theta^j}b^ib^j}{g_{ij}b^ib^j}.
\end{multline}
Let take a map:
$$\exp^{-1}_{\mbox o}: M^{n+1}\to T_{\mbox o}M^{n+1}=E^{n+1}.$$
The image of the sphere $S(R)$ is the sphere $\bar{S}(R)$ with the
center $\bar{\mbox{O}}~=~\exp^{-1}(\mbox{O})$ and radius $R$. The
image of the sphere $S(r)$ lies in a clo\-sed ball in  $E^{n+1}$
of radius $r$ with the center $\bar{P}=\exp^{-1}_{\mbox{o}}(P)$,
where $P$ is the center of the sphere $S(r)$.

\quad  Really, let consider triangles $OPX,\bar{O}\bar{P}\bar{X}$,
where

$X\in S(r)$, $\bar{X}=\exp^{-1}_{\mbox{o}}(X)$;
$OP=\bar{O}\bar{P}=R-r,\quad OX=\bar{O}\bar{X}=h$ and $\angle
POX=\angle\bar{P}\bar{O}\bar{X}$. From nonpositivity of the
sectional curvature $M^{n+1}$ and comparison theorem for triangles
if follows that $\bar{P}\bar{X}\leq PX$. In spherical system of
coordinates with pole $\bar{O}$ the metric $E^{n+1}$ has the form
$$ds^2=dt^2+G_{ij}(t,\theta)d\theta^id\theta^j.$$
The normal curvature of the image of the sphere $S(r)$ at the
point $\bar{Q}$ is equal
\begin{multline}\label{Bor20}
{\overline
k_n}=\dfrac{\left(\dfrac{\partial^2h(\theta^1,\ldots,\theta^n)}
{\partial\theta^i\partial\theta^j} -\dfrac{1}{2}\dfrac{\partial
G_{ij}}{\partial
t}\right)b^ib^j}{G_{ij}b^ib^j}=\\\frac{\dfrac{\partial^2h(\theta^1,\ldots,\theta^n)}
{\partial\theta^i\partial\theta^j}b^ib^j}{G_{ij}b^ib^j}+\dfrac{1}{R^2}.
\end{multline}
As the image $S(r)$ lies in a closed ball of radius $r$ with
center  $\bar{P}$, then
$$\bar{k}_n\geqslant\dfrac{1}{r}.$$
From \eqref{Bor20} it follows at the point $Q$
\begin{equation}\label{Bor21}
\dfrac{\partial^2h(\theta)}{\partial\theta^i\partial\theta^j}b^ib^j>0
\end{equation}
From \eqref{Bor19} and \eqref{Bor21} we obtain the statement of
the lemma1.
\end{proof}

\quad It follows from lemma that normal curvatures of the
horosphere less than normal curvatures of the tangent sphere which
lies inside  horoball, bounded by horosphere. Therefore at the
point $Q_0\in F^n$ normal curvatures $F^n$ satisfy an inequality:
$$k_n/F^n\geqslant k_n/S(r)>k_n/H^n,$$
where $H^n$ is the supporting tangent horosphere. But this
contradicts  the assumption that at any point $F^n$  there exists
the direction $a$ such that
$$k_n(a)/_{F^n}=k_n(a)/_{H^n}.$$
And the case $b)$ is impossible. The case $a)$ is possible only
for $\tau^*=\infty$, otherwise it is true arguments of the case
$b)$.

\quad We have proved that any tangent horosphere is globally
supporting.

\quad Let $P_1,P_2$ be different arbitrary points  $F^n$ and
tangent supporting horospheres  $H_1,H_2$ are different too. Then
$F^n$ belongs to intersection of horoballs bounded by horosphere
$H_1,H_2$. Intersection of horoballs is a compact bounded set if
the sectional curvature of Hadamard manifold
$$K_{\sigma}\leq-k_1^2<0.$$
Therefore  $\tau^*<\infty$, but it is impossible. Hence horosphere
$H_1$ and $H_2$ coincide and $F^n$ is a horosphere in Hadamard
manifold $M^{n+1}$.

\item[3)] Let introduce the horospherical system of coordinates
with base $F^n$ in the manifold $M^{n+1}$. The metric of the
ambient space has the form~\eqref{bor2}. For  $t=0$ we obtain the
hypersurface $F^n$. Principal curvatures of horosphere $t=const$
satisfy the inequalities  $k_2\geqslant \lambda_i\geqslant k_1$,
such that the sectional curvature of $M^{n+1}$ satisfies
inequality
$$-k_1^2\geqslant K\geqslant-k_2^2.$$
\quad By condition of the theorem the principal curvatures of
$F^n$ satisfy inequality $\lambda_i\geqslant k_2$ and we obtain
that $\lambda_i=k_2$ and horosphere $F^n$ is an umbilical
hypersurface.
 Principal curvature of equidistant horospheres $t=const$ satisfy
 the Riccati equation.
$$\frac{d\lambda}{dt}=\lambda^2+K_\sigma,$$
where $K_\sigma$ is the sectional curvature in the direction of
twodimensional plane span on the normal to horosphere and
corresponding principal direction. Since
$K_\sigma\geqslant-k_2^2$, that
$$\frac{d\lambda}{dt}\geqslant\lambda^2-k_2^2,\qquad \lambda(0)=\lambda_0.$$
Solving this inequality we obtain
$$\lambda\geqslant k_2\dfrac{(k_2+\lambda_0)e^{-2k_2t}-(k_2-\lambda_0)}
{(k_2+\lambda_0)e^{-2k_2t}+(k_2-\lambda_0)}$$ for
$\lambda_0=k_2,\quad \lambda\geqslant k_2$, from another side
$\lambda\leq k_2$. And we get $ \lambda=k_2$ for all values $t$.

\quad Therefore the coefficients of metric tensor $g_{ij}$ of the
ambient space $M^{n+1}$ satisfies the equations:
$$-\dfrac{1}{2}\dfrac{\partial g_{ij}}{\partial t}=k_2g_{ij}.$$
And  $g_{ij}(\theta,t)=g_{ij}(\theta,0)e^{-2k_2t}$. The metric
$M^{n+1}$  has the form
$$ds^2=dt^2+e^{-2k_2t}d\sigma^2,$$
where $d\sigma^2$ is the metric of the base horosphere $F^n$. Let
show that metric of $F^n$ is flat. Suppose that in some point of
$F^n$ on some twodimensional plane the sectional curvature
$\gamma_2\neq0$.  Then the sectional curvatures of the coordinates
horosphere $t=const$ in corresponding point and direction is equal
$\gamma_2e^{2k_2t}$. From Gauss formula we get that the sectional
curvature of the ambient space $M^{n+1}$ at the same direction is
equal
$$\gamma_2e^{2k_2t}-k_2^2,\quad -\infty\leq t<+\infty.$$ As the
sectional curvature $M^{n+1}$ satisfies the inequality
$$-k_1^2\geqslant K\geqslant -k_2^2,$$ that $\gamma_2=0$ and the manifold
$M^{n+1}$ is a space of constant curvature $-k_2^2$.
\end{enumerate}
\end{enumerate}
  \end{proof1}
  \begin{proof2}
  From the part I) of theorem 1 it follows that $F^n$ is  a
  compact convex hypersurface diffeomorphic to $S^n$. Analogical,
  to the proof of the theorem 3.1 \cite{9}  we obtain that
  every tangent sphere of radius $r_0$ is globally supporting and
  $F^n$ belongs to closed balls bounded of this spheres. It is
  possible two  cases:
  \begin{enumerate}
  \item[I).]
  There exist two different points $P_1, P_2\in F^n$ such that
  tangent spheres $S_1(r_0), S_2(r_0)$ at these points of radius
  $r_0$ don't coincide. Than $F^n$  lies in intersection of balls
  bounded of these spheres. In Hadamard manifold the intersection
  of different balls of radius $r_0$ belongs to the ball of radius
  less $r_0$.
  \item[II).]
  At all points $F^n$ the tangent sphere of radius $r_0$ is the
  same and $F^n$ coincides  with the sphere of radius $r_0$.
  Analogical to the proof of part II).3) of theorem1 we obtain that the ball
  bounded of this sphere  isometric to a ball of radius $r_0$ in
  Lobachevsky space of curvature $-k^2_2$.
  \end{enumerate}
  \end{proof2}
\begin{center}
\renewcommand{\refname}{\bf {References.}}
  
 \end{center}
\end{document}